\documentclass[11pt,twoside]{article}

\usepackage{a4wide}
\usepackage{amsfonts}
\usepackage{amssymb}
\usepackage{amsmath}
\usepackage{graphicx}

\newcommand{\ignore}[1]{}

\def\@begintheorem#1#2{\par\bgroup{\sc #1\ #2. }\it\ignorespaces}
\def\@opargbegintheorem#1#2#3{\par\bgroup{\sc #1\ #2\ (#3). } \it\ignorespaces}
\def\@endtheorem{\egroup}
\newtheorem{theorem}{Theorem}[section]
\newtheorem{corollary}[theorem]{Corollary}
\newtheorem{lemma}[theorem]{Lemma}

\newtheorem{example}[theorem]{Example}
\newtheorem{proposition}[theorem]{Proposition}
\newtheorem{definition}[theorem]{Definition}
\newcommand{\bt}[1]{\begin{theorem}\label{#1}}
\newcommand{\bc}[1]{\begin{corollary}\label{#1}}
\newcommand{\bl}[1]{\begin{lemma}\label{#1}}
\newcommand{\be}[1]{\begin{example}\label{#1}}
\newcommand{\bp}[1]{\begin{proposition}\label{#1}}
\newcommand{\ba}[1]{\begin{algorithm}\rm\label{#1}}
\newcommand{\bd}[1]{\begin{definition}\rm\label{#1}}{\normalsize }
\newcommand{\bpr}{\noindent {\em Proof. }}
\newcommand{\et}{\end{theorem}}
\newcommand{\ec}{\end{corollary}}
\newcommand{\el}{\end{lemma}}
\newcommand{\ee}{\end{example}}
\newcommand{\ep}{\end{proposition}}
\newcommand{\ed}{\end{definition}}
\newcommand{\epr}{{\ \vbox{\hrule\hbox{%
\vrule height1.3ex\hskip0.8ex\vrule}\hrule}}\\\par}

\def\R{\mathbb{R}}
\def\Z{\mathbb{Z}}
\def\xx{\overline x}
\def \conv {{\rm conv}}
\def \nk {\{0,1\}^n_k}
\def \nt {\{0,1\}^n_3}
\def \nkm {\left(\{0,1\}_k^n\right)^m}

\def \sign {{\rm sign}}

\begin{document}

\title{\bf Optimization over Degree Sequences}

\author{
Antoine Deza
\thanks{\small McMaster University, Hamilton.
Email: deza@mcmaster.ca}
\and
Asaf Levin
\thanks{\small Technion - Israel Institute of Technology, Haifa.
Email: levinas@ie.technion.ac.il.}
\and
Syed M. Meesum
\thanks{\small Institute of Mathematical Sciences, HBNI, Chennai.
Email: meesum@imsc.res.in}
\and
Shmuel Onn
\thanks{\small Technion - Israel Institute of Technology, Haifa.
Email: onn@ie.technion.ac.il}
}
\date{}

\maketitle

\begin{abstract}
We introduce and study the problem of optimizing arbitrary
functions over degree sequences of hypergraphs and
multihypergraphs. We show that over multihypergraphs the problem
can be solved in polynomial time. For hypergraphs, we show that
deciding if a given sequence is the degree sequence of a
$3$-hypergraph is NP-complete, thereby solving a 30 year long open
problem. This implies that optimization over hypergraphs is hard
even for simple concave functions. In contrast, we show that
for graphs, if the functions at vertices are the same, then the
problem is polynomial time solvable. We also provide positive
results for convex optimization over multihypergraphs and graphs
and exploit connections to degree sequence polytopes and threshold
graphs. We then elaborate on connections to the emerging theory of
shifted combinatorial optimization.

\end{abstract}

\section{Introduction}

The {\em degree sequence} of a (simple) graph $G=(V,E)$
with $V=[n]:=\{1,\dots,n\}$ and $m=|E|$ edges is the vector
$d=(d_1,\dots,d_n)$ with $d_i=|\{e\in E:i\in e\}|$ the degree of
vertex $i$ for all $i$.

Degree sequences have been studied by many authors including
the celebrated work of Erd\H{o}s and Gallai \cite{EG} from 1960
which effectively characterizes the degree sequences of graphs. See for example \cite{EKM}
and the references therein for some of the more recent work in this area.

In this article we are interested in the following discrete optimization problem.
Given $n,m$, and functions $f_i:\{0,1,\dots,m\}\rightarrow\R$ for $i=1,\dots,n$, find a graph
on $[n]$ with $m$ edges having degree sequence $d$ maximizing $\sum_{i=1}^n f_i(d_i)$.
The special case with $f_i(t)=t^2$ for all $i$, that is, finding a graph
maximizing the sum of degree squares, was solved previously, in \cite{PPS}.

More generally, we are interested in the problem over (uniform)
hypergraphs. In the sequel it will be convenient to use a
vector notation and so we make the following definitions.
A {\em $k$-hypergraph with $m$ edges on $[n]$} is a subset
$H\subseteq\nk:=\{x\in\{0,1\}^n:\|x\|_1=k\}$ with $|H|=m$ (keeping
in mind also the interpretation of $H$ as the set of supports of
its vectors). We also consider the problem over multihypergraphs.
A {\em $k$-multihypergraph with $m$ edges on $[n]$}
is a matrix $H\in\nkm$, that is, an $n\times m$ matrix with each
column $H^j\in\nk$ representing an edge (that is, each column is a
$0-1$ vector containing exactly $k$ ones), so that multiple
(identical) edges are allowed (but no loops). The {\em degree
sequence of $H$} is the vector $d=\sum H:=\sum\{x:x\in H\}$ for
hypergraphs and $d=\sum H:=\sum_{j=1}^m H^j$ for multihypergraphs.

We pose the following algorithmic problem.

\vskip.2cm\noindent{\bf Optimization over Degree Sequences.}
Given $k,n,m$ and functions $f_i:\{0,\dots,m\}\rightarrow\R$,
find a $k$-(multi)hypergraph $H$ with $m$ edges
whose degree sequence $d:=\sum H$ maximizes $\sum_{i=1}^n f_i(d_i)$.

\vskip.2cm
The case of linear functions will be discussed in
Section 3. For arbitrary functions and multihypergraphs
we solve the problem completely in Section 4.

\bt{thm_multihypergraph}
The general optimization problem over degree sequences of multihypergraphs can be
solved in polynomial time for any $k,n,m$ and any univariate functions $f_1,\dots,f_n$.
\et

For hypergraphs the problem is much harder. On the positive side we show in Section 5 the
following theorem, broadly extending the result for sum of degree squares of graphs.

\bt{thm_graph} For $k=2$, that is, graphs, the optimization
problem over degree sequences can be solved in $O(n^5m^2)$ time
for any $n,m$ and any identical univariate functions
$f_1=\cdots=f_n$.
\et
When the identical functions are {\em convex} we establish better time complexity
in Theorem \ref{thm_convex_graph}.

\vskip.2cm
In order to obtain our result on the negative side,
we consider also the decision problem: given $k$ and $d\in\Z_+^n$,
is $d$ the degree sequence of some hypergraph $H\subseteq\nk$?
For $k=1$ it is trivial as $d$ is a degree sequence if and only if $d\in\{0,1\}^n$.
For $k=2$ it is solved by the aforementioned theorem of Erd\H{o}s and Gallai \cite{EG}
which implies that $d$ is a degree sequence of a graph
if and only if $\sum d_i$ is even and, permuting $d$ so that $d_1\geq\dots\geq d_n$,
the inequalities $\sum_{i=1}^j d_i-\sum_{i=l+1}^n d_i\leq j(l-1)$ hold for $1\leq j\leq l\leq n$,
yielding a polynomial time algorithm.

\vskip.2cm
For $k=3$ it was raised 30 years ago by Colbourn, Kocay and Stinson \cite{CKS} (Problem 3.1)
and remained open to date. We solve it in Section 2,
implying it is unlikely that degree sequences of $3$-hypergraphs
could be effectively characterized, in the following theorem.
\bt{complete}
It is NP-complete to decide if $d\in\Z_+^n$ is the degree sequence of a $3$-hypergraph.
\et

This leads at once to the following negative statement by presenting an optimization problem
of simple concave functions over degree sequences of $3$-hypergraphs whose optimal
value is zero if and only if a given $d\in\Z_+^n$ is a degree sequence of some $3$-hypergraph.

\bc{hard_multihypergraph}
For $k=3$ the optimization problem over degree sequences of $3$-hypergraphs is NP-hard even
for concave functions of the form $f_i(t)=-(t-d_i)^2$ with $d_i\in\Z_+$ for each $i$.
\ec

Next, in Section 6, we discuss optimization of convex functions over degree sequences.
In fact, our results hold for the more general problem of maximizing any
convex function $f:\{0,1,\dots,m\}^n\rightarrow\R$ which is not necessarily separable
as considered above. For this we discuss the {\em degree sequence polytopes}
studied in \cite{KR,Kor,Liu,MS,PS} and references therein, introduce and study degree
sequence polytopes of hypergraphs with prescribed number of edges, and show that for
$k=2$ their vertices correspond to suitable {\em threshold graphs} \cite{MP}.

Finally, in Section 7, we illustrate that optimization over degree
sequences can be viewed within the framework of {\em shifted
combinatorial optimization} recently introduced and investigated
in the series of papers \cite{GHKO,KOS,KLMO,LO}, and contributes
to this emerging new theory.

\section{The complexity of deciding hypergraph degree sequences }

Here we consider the complexity of deciding the existence of a hypergraph with a given degree
sequence. We prove the following theorem solving a problem raised 30 years ago by Colbourn,
Kocay and Stinson \cite{CKS}. The proof is partially inspired by an argument from \cite{Liu}.

\vskip.2cm
\noindent{\bf Theorem \ref{complete}\ }
It is NP-complete to decide if $d\in\Z_+^n$ is the degree sequence of a $3$-hypergraph.

\vskip.2cm
\bpr
The problem is in NP since if $d$ is a degree sequences then a hypergraph
$H\subseteq\nt$ of cardinality $|H|\leq{n\choose3}=O(n^3)$
can be exhibited and $d=\sum H$ verified in polynomial time.

Now consider the following so-called {\em $3$-partition} problem:
given $a\in\Z_+^n$ and $b\in\Z_+$, decide if there is an $H\subseteq\nt$
such that $ax:=\sum_{i=1}^n a_ix_i=b$ for all $x\in H$
and $\sum H={\bf 1}$ where ${\bf 1}$ is the all-ones vector.
It is well known to be NP-complete \cite{GJ} and we reduce it to ours.
(In \cite{GJ} the problem is given in an equivalent form in terms of sets
rather than vectors, and even the special case with $n=3m$ for some $m$,
$a{\bf 1}=mb$ and ${1\over4}b<a_i<{1\over2}b$ for all $i$ is NP-complete.)

Given such $a$ and $b$ define $w\in\Z^n$ by $w_i:=3a_i-b$ for all $i$. Then for any $x\in\nt$
we have $wx=3ax-b\sum_{i=1}^n x_i=3(ax-b)$ and so $wx=0$ if and only if $ax=b$.
So $H$ satisfies $ax=b$ for all $x\in H$ and $\sum H={\bf 1}$ if and only if $wx=0$
for all $x\in H$ and $\sum H={\bf 1}$. For this to hold we must have
$w{\bf 1}=w\sum H=\sum\{wx:x\in H\}=0$. So if $w{\bf 1}\neq 0$ then there is no solution to the
$3$-partition problem and we can define $d$ to be a unit vector in $\R^n$ making sure there
is no $3$-hypergraph with degree sequence $d$ as well. So we may and do assume $w{\bf 1}=0$.

For $\sigma\in\{-,0,+\}$ define a hypergraph
$S_\sigma:=\{x\in\nt\,:\,\sign(wx)=\sigma\}$ so that these three hypergraphs
form a partition $S_-\uplus S_0\uplus S_+=\nt$ of the complete $3$-hypergraph.

Define $d:={\bf 1}+\sum S_+$ and observe that $w$, $S_+$, and $d$
can be computed in polynomial time. We claim that there is a
$3$-hypergraph $G$ with degree sequence $d$ if and only if there
is a $3$-hypergraph $H$ with $wx=0$ for all $x\in H$ and $\sum
H={\bf 1}$. So the $3$-partition problem reduces to deciding
degree sequences, showing the latter is NP-complete.

We now prove the claim. Suppose first $H$ satisfies $wx=0$ for all
$x\in H$ and $\sum H={\bf 1}$. Then $H\subseteq S_0$ so $H\cap
S_+=\emptyset$. Let $G:=H\uplus S_+$. Then $\sum G=\sum H+\sum
S_+={\bf 1}+\sum S_+=d$ so $G$ has degree sequence $d$.
Conversely, suppose $G$ has degree sequence $d$. Then, using $w{\bf 1}=0$,
$$w\sum G\ =\ \sum_{x\in G\cap S_-}wx+\sum_{x\in G\cap S_0}wx+\sum_{x\in G\cap S_+}wx
\ \leq\ \sum_{x\in S_+}wx\ =\ w{\bf 1}+w\sum S_+\ =\ wd$$
with equality if and only if $G\cap S_-=\emptyset$ and $G\cap S_+=S_+$,
since $\sign(wx)=\sigma$ for each $\sigma\in\{-,0,+\}$ and every $x\in S_\sigma$.
Since $\sum G=d$ we do have equality $w\sum G=wd$ above. Let $H:=G\cap S_0$. Then
$wx=0$ for all $x\in H$ and $\sum H=\sum G-\sum S_+={\bf 1}$ as claimed.
\epr

This result can be easily extended to $k$-hypergraphs for all
fixed $k\geq 3$ using a reduction from the problem for
$3$-hypergraphs into the problem for $k$-hypergraphs.
\begin{corollary}
Fix any $k\geq 3$. It is NP-complete to decide if $d\in\Z_+^n$ is
the degree sequence of a $k$-hypergraph.
\end{corollary}
\bpr
Given an input $d=(d_1,d_2,\ldots,d_n)$ to the problem for $3$-hypergraphs,
let $m={1\over3}\sum_{i=1}^n d_i$ and define an input
$d'=(d'_1,d'_2,\ldots d'_{n+k-3})$ to the problem for $k$-hypergraphs where
$d'_i=d_i$ for $1\leq i \leq n$ and $d'_i=m$ for $i>n$.
Then any $k$-hypergraph with degree sequence $d'$ has $m$ edges and each
must contain vertices $n+1,\ldots ,n+k-3$. It can then be verified
that $d$ is a degree sequence of a $3$-hypergraph
if and only if $d'$ is a degree sequence of a $k$-hypergraph.
\epr

\section{Linear functions}

Here we discuss as a warm-up to the following sections the
case of linear functions, that is, with $f_i(t)=w_it$ for each
$i$, where $w=(w_1,\dots,w_n)$ is a given profit vector. So the
degree sequence optimization problems are
$$\max\{w\sum H\ :\ H\in\nkm\}\ ,\quad\quad \max\{w\sum H\ :\ H\subseteq\nk\,,\ |H|=m\}\ .$$
\bp{linear}
The linear optimization problem over degree sequences of multihypergraphs and over degree sequences of
hypergraphs can be solved in polynomial time for all $n,k,m$.
\ep
\bpr
First, consider the case of multihypergraphs. Since $w\sum H=\sum_{j=1}^m wH^j$ and all columns $H^j$
can be the same, an optimal solution will be a matrix $H=[x,\dots,x]$ with all columns equal to
some $x\in\nk$ maximizing $wx$. Such an $x$ can be found by sorting $w$, that is,
if for a permutation $\pi$ of $[n]$ we have
$w_{\pi(1)}\geq\cdots\geq w_{\pi(n)}$, then we can take
$x:={\bf 1}_{\pi(1)}+\cdots+{\bf 1}_{\pi(k)}$
where ${\bf 1}_i$ denotes the standard $i$th unit vector in $\R^n$.

Next, consider the case of hypergraphs. Since $w\sum H=\sum\{wx:x\in H\}$, an optimal solution will be a
hypergraph $H$ consisting of the $m$ vectors $x\in\nk$ with largest values $wx$.
For $k$ fixed we can simply compute $wx$ for each of the $O(n^k)$ vectors in $\nk$
in polynomial time and pick the best $m$ vectors. For variable $k$, we can use the
algorithm of Lawler \cite{Lawler72} to find the $m$ vectors $x\in\nk$ with largest values $wx$
in time polynomial in $n,k,m$.
\epr

\section{Arbitrary functions over multihypergraphs}

The algorithmic results of our work are mostly obtained by the use
of dynamic programming.  We will use a model of dynamic
programming, where there is a finite set of {\em states} that is
partitioned into {\em stages} (the stages are indexed by $\{
1,2,\ldots ,n\}$). For each state, there is a finite set of {\em
actions} that can be performed. Given a state and an action, the
decision maker is moving to a state of a (strictly) larger stage,
and incurs a {\em reward} that is a function of the state and the
action.  We have an initial state, and a set of final states, and
the goal is to find a path maximizing the total reward starting in
the initial state and ending at a final state.  This optimization
problem is equivalent to finding a longest path in a
directed acyclic graph, a task known to be solvable with time
complexity that is linear in the number of arcs and nodes of that
directed graph.

For multihypergraphs there is a characterization of degree sequences and
a greedy procedure for constructing a hypergraph from its degree sequence.
The characterization and algorithm follow from results of \cite{Rys} on $0-1$
matrices and were also proved in \cite{Bil1,Bil2}. We record these facts
in the following proposition and provide a short proof for completeness.

\bp{multihypergraph_characterization} Vector $d\in\Z_+^n$ is a
degree sequence of $k$-multihypergraph $H$ with $m$ edges if and
only if $\sum_{i=1}^n d_i=km$ and $d_i\leq m$ for all $i$,
and $H$ is polynomial time realizable  from $d$.
\ep

\bpr The conditions on $d$ are clearly necessary. We prove by
induction on $m$ that given $d$ satisfying the conditions we can
construct a multihypergraph $H$ with $d=\sum H$. If $m=0$ this is
the empty multihypergraph. Suppose $m\geq1$. Permuting the $d_i$
we may assume $d_1\geq\cdots\geq d_n$. We claim $d_k\geq 1$ and
$d_{k+1}\leq m-1$. If $d_k=0$ then $km=\sum d_i\leq(k-1)m$ which
is impossible. If $d_{k+1}\geq m$ then $km=\sum d_i\geq(k+1)m$
which is again impossible. Define a vector $d'$ by $d'_i:=d_i-1$
for $i\leq k$ and $d'_i:=d_i$ for $i>k$. This $d'$ satisfies the
conditions with $m-1$ and by induction we can construct a
multihypergraph $[H^1,\dots,H^{m-1}]$ with degree sequence $d'$.
Then $H:=[H^1,\dots,H^{m-1},H^m]$ with $H^m:={\bf 1}_1+\cdots+{\bf
1}_k$ is the desired multihypergraph with degree sequence $d$. So
the induction follows and we are done. \epr

We can now prove our theorem on optimization over degree sequences of multihypergraphs.
\vskip.2cm
\noindent{\bf Theorem \ref{thm_multihypergraph}\ }
The general optimization problem over degree sequences of multihypergraphs can be
solved in polynomial time for any $k,n,m$ and any univariate functions $f_1,\dots,f_n$.

\vskip.2cm
\bpr
By Proposition \ref{multihypergraph_characterization},
we need to solve the following integer programming problem,
$$\max\{\sum_{i=1}^n f_i(d_i)\ :\
d\in\Z_+^n\,,\ \sum_{i=1}^n d_i=km\,,\ d_i\leq m\,,\ i=1,\dots,n\}\ ,$$
find an optimal $d$, and then use the algorithm of
Proposition \ref{multihypergraph_characterization} to find an $H$ with $d=\sum H$.

An optimal $d$ can be found using dynamic programming.
There are $n$ stages, where the decision at stage $i$ is $0\leq d_i\leq m$
with reward $f_i(d_i)$. The state at the end of stage $i$ is $0\leq s_i\leq km$ representing
the partial sum $s_i=\sum_{j=1}^i d_j$ starting with $s_0:=0$.
Let $f^*_i(s_i)$ be the maximum total reward of a path leading from the initial
state $s_0$ to state $s_i$ in stage $i$.
The recursively computable optimal value at state $s_i$ is given by
$$f^*_0(s_0)\ :=0\ ,\quad f^*_i(s_i)
\ :=\ \max\{f^*_{i-1}(s_{i-1})+f_i(d_i)\ :\ s_i=s_{i-1}+d_i\}\ ,\quad i=1,\dots,n\ .$$
The optimal value is $f^*_n(km)$ and the optimal decisions and the sequence of states on the
optimal path can be reconstructed backwards starting with $s^*_n:=km$ recursively by
$$(s^*_{i-1},d^*_i)\ \in\
\arg\max\{f^*_{i-1}(s_{i-1})+f_i(d_i)\ :\ s^*_i=s_{i-1}+d_i\}\
,\quad i=n,\dots,1\ .$$ Clearly this is doable in time polynomial
in $k,n,m$ for any univariate functions $f_1,\dots,f_n$. This last
claim follows since there are $O(n)$ stages and in each stage
there are $O(mk)$ states and the computation of $f^*_i$ for a
given state in stage $i$ takes $O(m)$ (for all $i$), and thus the
time complexity of this algorithm is $O(nkm^2)$. \epr

\section{Identical functions over graphs}

In this section we restrict attention to graphs. We use alternatively the interpretation
of a graph as $G\subseteq\{0,1\}^n_2$ and $G=([n],E)$.
We note the following characterization of degree sequences of graphs from \cite{Hak,Hav}
which leads to a greedy procedure for constructing a graph from its degree sequence.
A vector $d\in\Z_+^n$ is called {\em reducible} if, permuting it such that $d_1\geq\cdots\geq d_n$,
we have that $d_{d_1+1}\geq 1$. The {\em reduction} of $d$ is the vector $d'\in\Z_+^{n-1}$
defined by $d':=(d_2-1,\dots,d_{d_1+1}-1,d_{d_1+2},\dots,d_n)$.
We record the characterization and algorithm in the next proposition
and provide a short proof for completeness, see e.g. \cite{Wes} for more details.

\bp{graph_characterization} Vector $d\in\Z_+^n$ is a degree
sequence of graph $G$ if and only if $d$ is reducible with
reduction also a degree sequence. So $G$ is polynomial time recursively
realizable from $d$.
\ep
\bpr Assume $d_1\geq\cdots\geq d_n$. Suppose $d$ is reducible and
its reduction $d'$ is the degree sequence of a graph $G'$ on
$[n]\setminus\{1\}$. Then the graph $G$ obtained from $G'$ by
adding vertex $1$ and connecting it to vertices $2,\dots,d_1+1$
has degree sequence $d$.

Conversely, suppose $d$ is the degree sequences of $G=([n],E)$. We show it is also the degree
sequence of a graph with $1$ as a neighbor of $2,\dots,d_1+1$, which will show that $d$ is reducible
and its reduction is also a degree sequence. Call $i\neq j$ a {\em bad pair} if $\{1,i\}\in E$
and $\{1,j\}\not\in E$ but $d_i<d_j$. Then there must be a $k$ with  $\{i,k\}\not\in E$
and $\{j,k\}\in E$. Then the graph obtained from $G$ by dropping edges
$\{1,i\}$, $\{j,k\}$ and adding $\{1,j\}$, $\{i,k\}$ has the same degree sequence but
fewer bad pairs. Repeating this procedure we arrive at a graph with degree sequence $d$
and no bad pairs which implies that it has $1$ as a neighbor of $2,\dots,d_1+1$.
\epr

We can now prove our theorem on optimization over degree sequences of graphs.
\vskip.2cm
\noindent{\bf Theorem \ref{thm_graph}\ }
For $k=2$, that is, graphs, the optimization problem over degree sequences can be solved in
$O(n^5m^2)$ time for any $n,m$ and any identical univariate functions $f_1=\cdots=f_n$.

\vskip.2cm
\bpr
Since the functions $f_i$ are assumed identical, we may optimize over {\em sorted}
degree sequences $d_1\geq\cdots\geq d_n$, and in fact, this is the reason our approach
needs this assumption.  For simplicity of notation we denote by $f$ the common function $f_1=\cdots =f_n=f$.
By Proposition \ref{graph_characterization} and by the
original version of the Erd\H{o}s-Gallai characterization of sorted degree sequences,
we need to solve the following nonlinear integer optimization problem,
\begin{eqnarray*}
\max && \sum_{i=1}^n f(d_i) \\
s.t. && d\in\Z_+^n\,,\ \ d_1\geq\cdots\geq d_n\,,\ \ \sum d_i=2m \\
     && \sum_{i=1}^j d_i-j(j-1)\leq\sum_{i=j+1}^n\min\{j,d_i\}\,,\ \ j=1,\dots,n\ ,
\end{eqnarray*}
find an optimal $d$, and then use the algorithm of Proposition
\ref{graph_characterization} to find a graph $G$ with $d=\sum G$.

Finding an optimal $d$ can be done again by dynamic programming.
This time in every computational path there are at most $n$ stages,
where the decision in stage $i$ is the value of $d_i$ as well as the number $\alpha_i$
of vertices with degree exactly $i$ (of indices $\beta_i+1,\ldots \beta_i+\alpha_i$).
The reward of this decision is $f(d_i)+\alpha_i\cdot f(i)$.
A {\em final state} in the dynamic programming table is a state corresponding
to a value of $i$ that equals $\max \{ j : d_j >j\}$ (in this case $i$ is the largest index of a vertex with degree larger than $i$ and thus $\beta_i \geq i$ while $i+1$ has degree at most $i+1$, and thus by definition $\beta_{i+1} \leq i$ and by the monotonicity of the sequence of $\beta$'s we conclude that $\beta_i=i$ and furthermore a final state means that
the total degrees of all vertices is $2m$).
Our goal is to find a computational path (corresponding to a path of decisions)
that maximizes the total reward of the decisions along the path starting at the
initial state and ending at a final state.

The definition of the states as well as the transition function will enforce both the
monotonicity constraints and the Erd\H{o}s-Gallai  constraints.
Thus, the states of the dynamic programming table correspond to 5-tuples $(i,p_i,d_{i},\beta_i,s_i)$.
The meaning of reaching this state is that so far we decided upon the
degrees of nodes $1,2,...,i$ all of them are at least $d_{i}$
(as we enforce the monotonicity constraints by induction on the length of the path)
and their sum is $p_i$, and we also computed the degrees of nodes $\beta_i+1,\beta_i+2,\ldots,n$
which are the vertices of degrees at most $i$, and their sum is $s_i$.

When we are at state $(i,p_i,d_{i},\beta_i,s_i)$ such that $i<
\beta_i$ we should decide the value of $d_{i+1}$ and the value of
$\beta_{i+1}$. As mentioned above a feasible value of the pair
$(d_{i+1}, \beta_{i+1})$ has a reward of $f(d_{i+1})+(\beta_i -
\beta_{i+1})\cdot f(i+1)$. If such a decision is feasible then it
will results in a transition to the state $(i+1,
p_i+d_{i+1},d_{i+1},\beta_{i+1},s_i+(i+1)\cdot(\beta_i -
\beta_{i+1}))$.

It remains to define the set of pairs $(d_{i+1}, \beta_{i+1})$ of non-negative integers
for which the decision is feasible.  A pair $(d_{i+1}, \beta_{i+1})$ is {\em feasible}
for a given state $(i,p_i,d_{i},\beta_i,s_i)$ if all the following conditions hold:
\begin{enumerate}

\item $\beta_{i+1} \geq i+1$ - this condition is required as we determine the degrees of both a prefix
and a suffix of the vertices and we need this definition to be well-defined.  Using this condition we enforce that the degree of a vertex will
be defined at most once either as part of the prefix or as a part of the suffix.

\item $i+1\leq d_{i+1} \leq d_i$ - this constraint enforces the monotonicity conditions
over the prefix of the first $i+1$ vertices (note that $d_i$ is defined in the definition of the state).

\item  $p_i+d_{i+1}+ s_i+(i+1)\cdot (\beta_i - \beta_{i+1}) \leq
m$ - this condition will enforce that we can choose at most $m$
edges.

\item The $(i+1)$-th Erd\H{o}s-Gallai constraint that can be stated as follows
$$p_i+d_{i+1} - (i+1)\cdot i \leq s_{i}+ (i+1)\cdot (\beta_i -
\beta_{i+1}) + (\beta_{i+1}-i-1)\cdot (i+1) .$$

\end{enumerate}

The {\em initial state} of the dynamic programming is $(0,0,n,n,0)$ and we would like to find
a path of feasible decisions with maximum total reward leading from the initial state to a final
state that is a state for which both $i=\beta_i$ and $p_i+s_i=m$.

Observe that if the decision is feasible and leads to a state that has a path $P$ to a final
state then the decisions that will correspond to the path $P$ will not change the fact that
the degree sequence satisfies the $(i+1)$-th Erd\H{o}s-Gallai constraint.
In order to establish this property we examine the change in the two sides of the  $(i+1)$-th Erd\H{o}s-Gallai constraint when we determine the  degrees in the next stages.  First, consider the left hand side and note that the left hand side of this inequality is not changed as all values that appear in the left hand side were determined in the previous stages. Consider the right hand side of this inequality and we argue that it also does not change.  This claim holds as all the degrees that we select in the following
stages will be between $i+2$ and $d_{i+1}$ so the minimum on the right hand side for the indices
between $i+2$ and $\beta_{i+1}$ will be $i+1$, and this is the value we used in the last condition we checked for
the feasibility of the decision.

Next we note that the number of states of the dynamic programming table is $O(n^3m^2)$ as
there are $O(n)$ options for $i$, $O(m)$ options for $p_i$, $O(n)$ options for $d_i$, $O(n)$ options
for $\beta_i$, and $O(m)$ options for $s_i$. Furthermore for each state there are $O(n^2)$ feasible
decisions, and the time for computing the maximum total reward of a subpath that starts at the given
state and ends at some final state is $O(n^2)$ (since given a value of $\beta_{i+1}$ we can compute
the upper bound on $d_{i+1}$ resulting from all our constraints in time $O(n)$).
Thus, the time complexity of the algorithm is $O(n^5m^2)$.
\epr

\section{Convex functions and degree sequence polytopes}

Here we consider optimization of convex functions over degree sequences.
In fact, our results hold for the more general problem of maximizing any
convex function $f:\{0,1,\dots,m\}^n\rightarrow\R$ which is not necessarily separable
of the form $f(d_1,\dots,d_n)=\sum_{i=1}^n f_i(d_i)$ considered before.

First we consider the case of multihypergraphs.
\bt{convex}
If $f$ is convex then there exists a multihypergraph $H=[x,\dots,x]$, having $m$ identical
edges $x$, which maximizes $f(\sum H)$, where $x$ maximizes $f(mx)$ over $\nk$. So for any
fixed $k$ the optimization problem over degree sequences can be solved in polynomial time.
On the other hand, for $k$ variable the problem may need exponential time even for $m=1$.
\et
\bpr
Assume $f$ is convex and let $\hat H$ be an optimal multihypergraph. Let $x^1,\dots,x^r$ in
$\nk$ be the distinct edges of $\hat H$. Define an $n\times r$ matrix $M:=[x^1,\dots,x^r]$
and a function $g:\Z^r\rightarrow\R$ by $g(y):=f(My)$, which is convex since $f$ is.
Consider the integer simplex $Y:=\{y\in\Z_+^r\,:\,y_1+\cdots+y_r=m\}$ and the auxiliary
problem $\max\{g(y)\,:\,y\in Y\}$. For each $y\in Y$ let
$H(y):=[x^1,\dots,x^1,\dots,x^r,\dots,x^r]$ be the multihypergraph
consisting of $y_i$ copies of $x^i$. Permuting the columns of $\hat H$
we may assume that $\hat H=H(\hat y)$ for some $\hat y\in Y$.

Now, $g$ is convex, so the auxiliary problem has an optimal solution which is a vertex
of $Y$, namely, a multiple $\tilde y=m{\bf 1}_i$ of a unit vector in $\Z^r$.
It then follows that $\tilde H:=H(\tilde y)=[x^i,\dots,x^i]$ is the desired
optimal multihypergraph having a degree sequence which maximizes $f$, since
$$f\left(\sum_{j=1}^m\tilde H^j\right)\ =\ f(M\tilde y)\ =\ g(\tilde y)
\ \geq\ g(\hat y)\ =\ f(M\hat y)\ =\ f\left(\sum_{j=1}^m\hat H^j\right)\ .$$
For $k$ fixed we can clearly find an $x$ attaining $\max\{f(mx)\,:\,x\in\nk\}$ in polynomial time.

Now consider the situation of $k$ variable part of the input. Let $m:=1$. Then any
function $f:\{0,1\}^n\rightarrow\R$ is convex so in order to solve $\max\{f(x):x\in\nk\}$
one needs to check the value of $f$ on each of the ${n\choose k}$ points in $\nk$
and for $k=\lfloor{n\over2}\rfloor$ this requires exponential time.
\epr

\vskip.2cm
We continue with hypergraphs. We need to discuss the class of {\em degree sequence polytopes}.
The classical degree sequence polytope is the convex hull $D^n_k:=\conv\{\sum H: H\subseteq\nk\}$
of degree sequences of $k$-hypergraphs on $[n]$  with unrestricted number of edges.
For $k=2$, that is, graphs, these polytopes have been extensively studied, see \cite{PS}
and the references therein. The Erd\H{o}s-Gallai theorem implies that $D^n_2$ is the set
of points $d\in\R^n$ satisfying the system
$$\sum_{i\in S} d_i-\sum_{i\in T} d_i\ \leq\ |S|(n-1-|T|)
\ ,\quad S,T\subseteq[n]\,,\ S\cap T=\emptyset\ ,$$
and the vertices of $D^n_2$ were characterized in \cite{Kor} as precisely the degree sequences
of threshold graphs. More recently, the polytopes $D^n_k$ for $k\geq 3$ were
studied in \cite{KR,Liu,MS}, but neither a complete inequality description nor a
complete characterization of vertices is known.

We now go back to the situation when the number of edges $m$ is prescribed. We define the
{\em degree sequence polytope} $D^{n,m}_k$ as the convex hull
$D^{n,m}_k:=\conv\{\sum H: H\subseteq\nk\,,\ |H|=m\}$ of degree sequence of $k$-hypergraphs
with $m$ edges over $[n]$. A degree sequence is called {\em extremal} if it is a
vertex of this polytope. We note that when maximizing a convex function over degree
sequences there will always be an optimal solution which is extremal.

Again we restrict attention to graphs and use alternatively the interpretation of a graph
as $G\subseteq\{0,1\}^n_2$ and $G=([n],E)$. The graph $G$ is called a {\em threshold graph}
if for some permutation $\pi$ of $[n]$, each vertex $\pi(i)$ is connected to either all
or none of the vertices $\pi(j)$, $j<i$. It is well known that a graph $G$ with degree
sequence $d=\sum G$ is threshold if and only if $N(i)\subseteq N[j]$ whenever $d_i\leq d_j$,
where $N(i):=\{j\,:\,\{i,j\}\in E\}$ and $N[i]:=N(i)\uplus\{i\}$ are the open and closed
neighborhoods of $i$ respectively, see \cite{MP}. We now characterize the extremal degree sequences
of graphs with prescribed number of edges in analogy with the result of \cite{Kor}.

\bt{threshold}
The vertices of $D^{n,m}_2$ are the degree sequences of threshold graphs with $m$ edges.
So for any convex function $f:\{0,1,\dots,m\}^n\rightarrow\R$ there is a threshold graph
$G$ maximizing $f(\sum G)$ and providing an optimal solution to the optimization problem
over degree sequences.
\et
\bpr
Suppose first that $d$ is the degree sequence of a threshold graph with $m$ edges. Then,
by the result of \cite{Kor}, it is a vertex of $D^n_2$. Since $d\in D^{n,m}_2\subseteq D^n_2$,
it is also a vertex of $D^{n,m}_2$.

Conversely, consider any vertex $d$ of $D^{n,m}_2$. Let $w\in\R^n$ be such that $wx$ is
uniquely maximized over $D^{n,m}_2$ at $d$. Perturbing $w$ if necessary we may assume
that $w_1,\dots,w_n$ are distinct. Let $G=([n],E)$ be a graph with $d=\sum G$. We need to show
that $N(i)\subseteq N[j]$ whenever $d_i\leq d_j$, which will imply that $G$ is a threshold graph.
Suppose on the contrary that for some $i\neq j$ we have $d_i\leq d_j$ but there is some
$k\in N(i)\setminus N[j]$.

Suppose first $d_i=d_j$ and $w_i>w_j$. Then there must be also some
$l\in N(j)\setminus N[i]$. Then the graph $G'$ obtained from $G$ by adding edge
$\{i,l\}$ and dropping edge $\{j,l\}$ has $m$ edges and $w\sum G'=w\sum G+w_i-w_j>wd$
which is impossible.

Next, consider the case where $d_i \leq d_j$ and $w_j>w_i$.
Then the graph $G''$ obtained from $G$
by adding $\{j,k\}$ and dropping $\{i,k\}$ has $m$ edges and
$w\sum G''=w\sum G+w_j-w_i>wd$ which is impossible again.

Thus, it remains to consider the case where $d_i<d_j$ and $w_i>w_j$.
Let $d'\in D^{n,m}_2$ be the degree sequence obtained from $d$ by setting
$d'_i:=d_j$, $d'_j:=d_i$, and $d'_t:=d_t$ for $t\neq i,j$.
Then $wd'-wd=(w_i-w_j)(d_j-d_i)>0$ which is a contradiction.

The second statement of the theorem now also follows, since any convex function attains its
maximum over a polytope at a vertex, which is the degree sequence of a threshold graph.
\epr

Two remarks are in order here. First, we note that if the function is not convex then,
even if it is separable, it may be that no threshold graph is a maximizer. To see this,
let $m=n$ and define univariate functions $f_1=\cdots=f_n$ by $f_i(2)=1$ and $f_i(t)=0$
for all $t\neq 2$. Then any graph $G$ with $m$ edges whose degree sequence $d=\sum G$
maximizes $\sum_{i=1}^n f_i(d_i)$ must have degree sequence satisfying $d_1=\cdots=d_n=2$
and so must be a (vertex) disjoint union of circuits. For all $n\geq 4$ any such graph is
not threshold, since by definition any threshold graph contains either an isolated vertex
(of degree 0) or a dominating vertex (of degree $n-1\geq 3$).

Second, given a convex function $f$, we do not know how to efficiently find a vertex $d$
of $D^{n,m}_2$ which maximizes $f$. Note that if we could find such a vertex $d$ then,
since the above Theorem \ref{threshold} guarantees that $d$ must be the degree sequence of
some threshold graph, we could also find an optimal threshold graph $G$ with $\sum G=d$.
Indeed, there must be some $i$ with either $d_i=0$ or $d_i=n-1$. We then define
$V':=[n]\setminus\{i\}$. If $d_i=0$ then we recursively find a threshold graph $G'$ on $V'$
with degree sequence $d'$ defined by $d'_i:=d_i$ for all $i\in V'$, and add an isolated
vertex $i$. If $d_i=n-1$ we recursively find a threshold graph $G'$ on $V'$ with degree
sequence $d'$ defined by $d'_i:=d_i-1$ for all $i\in V'$, and add a dominating vertex $i$.

We conclude this section by showing that we can use Theorem \ref{threshold}
to obtain a dynamic programming algorithm for the optimization problem over degree
sequences for any identical univariate functions $f_1=\cdots=f_n=f$ which are convex.
Observe that the time complexity that is established in the next theorem is
significantly lower than the one established for the general case of (not necessarily convex)
identical univariate functions in Theorem \ref{thm_graph}.

\bt{thm_convex_graph} For $k=2$, that is, graphs, optimization
over degree sequences can be done in $O(n^2m)$ time for
any $n,m$ and any identical univariate convex functions
$f_1=\cdots=f_n$.
\et
\bpr
We use the fact that the functions $f_1=\cdots=f_n=f$ are identical to
assume that an optimal solution of our problem (which can be assumed to be a
threshold graph by Theorem \ref{threshold}) satisfies that every vertex $i$
is either adjacent to all vertices $j$ such that
$j>i$ (in which case we say that $i$ is {\em dominating}) or it is
not adjacent to any of those vertices (in which case we say that $i$
is {\em isolating}).  Thus, the terms dominating and isolating refer
to the induced subgraph over the vertices with indices at least $i$.

Once again we use dynamic programming to find an optimal threshold graph for our problem.
There are $n$ stages, and in stage $i$ we decide if vertex $i$ is dominating or isolating.
A state in this dynamic program is a triple consisting of $(i,e_i,\delta_i)$ where $e_i$
is the total number of edges adjacent to at least one vertex in $\{ 1,2,\ldots ,i\}$, and $\delta_i$
is the number of dominating vertices among $1,2,\ldots ,i$. The decision that $i+1$ is dominating
means that its degree will be $\delta_i+n-(i+1)$ as $i+1$ will be adjacent to all dominating
vertices in $1,\ldots,i$ and to all vertices in $i+2,\ldots ,n$, and thus the reward of such a
decision will be $f(\delta_i+n-(i+1))$ and in this case we move to state $(i+1,\delta_i+1,e_i+n-(i+1))$.
The decision that $i+1$ is isolating means that its degree will be $\delta_i$ and thus the reward
of such a decision will be $f(\delta_i)$ and we will move to the state $(i+1,\delta_i,e_i)$.

These decisions are feasible only when the third component is at most $m$
(and the second component is at most $n$), and otherwise we have only one option of deciding
that the vertex is isolated.

In the resulting dynamic programming we look for a maximum total reward path that leads from
the initial state of $(0,0,0)$ to any of the final states defined as $(n,\delta,m)$
(choosing the value of $\delta$ that maximizes the total reward of the path).

The time complexity of this algorithm is $O(n^2m)$ as there are $O(n^2m)$ states
(since there are $O(n)$ options for $i$, $O(n)$ options for $\delta_i$, and $O(m)$ options
for $e_i$) and the amount of work for each state is $O(1)$.
\epr

We conclude this section with the following proposition on the complexity of deciding
membership in degree sequence polytopes and deciding being extremal for these polytopes.

\bp{vertices}
It is polynomial time decidable, given $n,k,m$ and rational $d\in\R^n$, whether or not $d$
lies in the degree sequence polytope $D^{n,m}_k$ and whether or not $d$ is a vertex of $D^{n,m}_k$.
\ep
\bpr
By the polynomial time equivalence of optimization and separation via the ellipsoid method, if
we can do linear optimization over a rational polytope $P\subset\R^n$ in polynomial time then we can also
decide in polynomial time whether or not a given $d\in\R^n$ lies in $P$ and whether or not $d$
is a vertex of $P$, see \cite{GLS} for more details. Since by Proposition \ref{linear} we can do
linear optimization over $D^{n,m}_k$ for any $n,k,m$ in polynomial time, the proposition follows.
\epr

This proposition indicates that for any $n,k,m$, effective characterizations
of inequalities defining $D^{n,m}_k$ and vertices of $D^{n,m}_k$ are plausible;
these remain challenging research problems.

\section{Shifted combinatorial optimization}

Optimization over degree sequences can be viewed as a special case of the
broad framework of {\em shifted combinatorial optimization} recently
introduced and investigated in \cite{GHKO,KOS,KLMO,LO}.

Standard combinatorial optimization is the following extensively studied problem (see \cite{Sch}
for a detailed account of the literature and bibliography of thousands of articles on this).

\vskip.2cm\noindent{\bf (Standard) Combinatorial Optimization.}
Given $S\subseteq\{0,1\}^n$ and $w\in\R^n$, solve
\begin{equation}\label{standard}
\max\{ws\ :\ s\in S\}\ .
\end{equation}
The complexity of this problem depends on the type and presentation of the defining set $S$.

Shifted combinatorial optimization is a broad nonlinear extension of this problem, which
involves the choice of several feasible solutions from $S$ at a time, defined as follows.

For a set $S\subseteq\R^n$ let $S^m$ be the set of $n\times m$ matrices $x$ having each
column $x^j$ in $S$. Call matrices $x,y\in\R^{n\times m}$ {\em equivalent} and write $x\sim y$
if each row of $x$ is a permutation of the corresponding row of $y$. The {\em shift} of
$x\in\R^{n\times m}$ is the unique matrix $\xx\in\R^{n\times m}$ satisfying $\xx\sim x$ and
$\xx^1\geq \cdots\geq \xx^m$. We can then define the following broad optimization framework.

\vskip.2cm\noindent{\bf Shifted Combinatorial Optimization.}
Given $S\subseteq\{0,1\}^n$ and $c\in\R^{n\times m}$, solve
\begin{equation}\label{shift}
\max\{c\xx\ :\ x\in S^m\}\ .
\end{equation}

\vskip.2cm
This framework has a very broad expressive power and its polynomial time solvability was so
far established in the following situations:
\begin{itemize}
\item
when $S$ is presented by totally unimodular inequalities, in particular when $S$ is the set
of source-sink dipaths in a digraph or matchings in a bipartite graph \cite{KOS};
\item
when $S$ is the set of independent sets in a matroid, in particular
spanning trees in a graph, or the intersection of two
so-called strongly-base-orderable matroids \cite{LO};
\item
when $S$ is any property defined by any bounded monadic-second-order-logic
formula over any graph of bounded tree-width \cite{GHKO}.
\end{itemize}
Finally, a study of approximation algorithms for this and the closely related
{\em separable multichoice optimization}
framework over monotone systems was very recently taken in \cite{KLMO}.

\vskip.2cm
We claim that the optimization problem over degree sequences in multihypergraphs can
be formulated as the special shifted combinatorial optimization problem with $S=\nk$. Then
multihypergraphs are matrices $x\in S^m=\nkm$ with degree sequence $\sum x=\sum_{j=1}^m x^j$.

Now, given functions $f_i:\{0,1,\dots,m\}\rightarrow\R$, define $c\in\R^{n\times m}$
by $c_{i,j}:=f_i(j)-f_i(j-1)$ for all $i$ and $j$. Then for every $x\in S^m$ we have
$c\xx=\sum_{i=1}^n f_i(d_i)-\sum_{i=1}^nf_i(0)$ where $d=\sum x$. Therefore the
multihypergraph $x$ maximizes $c\xx$ if and only if it maximizes $\sum_{i=1}^n f_i(d_i)$.

\vskip.2cm
Next, for $S\subseteq\{0,1\}^n$ and positive integers $m,u$ let $S_u^m$ be the set of matrices
$x$ in $S^m$ such that for each $s\in S$ there are at most $u$ columns of $x$ which are
equal to $s$. The {\em bounded} shifted combinatorial optimization problem is then
$\max\{c\xx\,:\,x\in S_u^m\}$. Then, defining $c$ from the functions $f_i$ as above,
the degree sequence problem for hypergraphs can be formulated as a bounded shifted
combinatorial optimization problem with $S=\nk$ as before and $u=1$.

Thus, our results in this paper for hypergraphs can be regarded as a first step
in the development of a theory of bounded shifted combinatorial optimization.

\section*{Acknowledgment}

A. Deza was partially supported by the Natural Sciences and Engineering Research
Council of Canada Discovery Grant programs (RGPIN-2015-06163).
A. Levin was partially supported by a Grant from
GIF - the German-Israeli Foundation for Scientific Research and Development
(grant number  I-1366-407.6/2016).
S.M. Meesum was partially supported by a grant at the Technion.
S. Onn was partially supported by the Dresner Chair at the Technion.


\begin{thebibliography}{}

\bibitem{Bil1}
Billington, D.: Degree multisets of hypergraphs.
Ph.D. Thesis, University of Melbourne (1982)

\bibitem{Bil2}
Billington, D.: Conditions for degree sequences to be realisable by $3$-uniform hypergraphs.
Journal of Combinatorial Mathematics and Combinatorial Computing 3:71--91 (1988)

\bibitem{CKS}
Colbourn, C.J., W.L. Kocay, W.L., Stinson, D.R.:
Some NP-complete problmes for hypergraph degree sequences.
Discrete Applied Mathematics 14:239--254 (1986)

\bibitem{EG}
Erd\H{o}s, P., Gallai, T.:
Graphs with prescribed degrees of vertices (in Hungarian).
Matematikai Lopak 11:264--274 (1960)

\bibitem{EKM}
Erd\H{o}s, P., Kir\'aly, Z., Mikl\'os, I.:
On the swap-distances of different realizations of a graphical degree sequence.
Combinatorics, Probability and Computing 22:366--383 (2013)

\bibitem{GHKO}
Gajarsk\'y J., Hlin\v{e}n\'y P., Kouteck\'y M., Onn S.:
Parameterized shifted combinatorial optimization.
Proceedings of the 23rd Annual International Computing and Combinatorics Conference,
Lecture Notes in Computer Science 10392:224--236 (2017)

\bibitem{GJ}
Garey, M.R., Johnson, D.S.:
Computers and Intractability. Freeman (1979)

\bibitem{GLS}
Gr\"otschel, M., Lov\'asz, L., Schrijver, A.:
Geometric Algorithms and Combinatorial Optimization. Springer (1993)

\bibitem{Hak}
Hakimi, S.L.:
On realizability of a set of integers as degrees of the vertices of a linear graph.
SIAM Journal 10:496--506 (1962)

\bibitem{Hav}
Havel, V.:
A remark on the existence of finite graphs.
\v{C}asopis pro P\v{e}stov\'an\'i Matematiky 80:477--480 (1955)

\bibitem{KOS}
Kaibel V., Onn S., Sarrabezolles P.:
The unimodular intersection problem.
Operations Research Letters 43:592--594 (2015)

\bibitem{KR}
Klivans, C., Reiner, V.:
Shifted set families, degree sequences, and plethysm.
The Electronic Journal of Combinatorics 15:\#R14 (2008)

\bibitem{Kor}
Koren, M.:
Extreme degree sequences of simple graphs.
Journal of Combinatorial Theory Series B 15:213--224 (1973)

\bibitem{KLMO}
Kouteck\'y M., Levin, A., Meesum S.M., Onn S.:
Approximate separable multichoice optimization over monotone systems.
Submitted.

\bibitem{Lawler72}
Lawler, E.L.:
A procedure for computing the $k$ best solutions to discrete optimization
problems and its application to the shortest path problem.
Management Science 18:401--405 (1972)

\bibitem{LO}
Levin A., Onn S.:
Shifted matroid optimization.
Operations Research Letters 44:535--539 (2016)

\bibitem{Liu}
Liu, R.I.:
Nonconvexity of the set of hypergraph degree sequences.
The Electronic Journal of Combinatorics 20:\#P21 (2013)

\bibitem{MP}
Mahadev, N.V.R., Peled, U.N.:
Threshold Graphs and Related Topics.
Annals of Discrete Mathematics 56, North-Holland (1995)

\bibitem{MS}
Murthy, N.L.B., Srinivasan, M.K.:
The polytope of degree sequences of hypergraphs.
Linear Algebra and its Applications 350:147--170 (2002)

\bibitem{PPS}
Peled, U.N., Petreschi, R., Sterbini, A.:
$(n,e)$-graphs with maximum sum of squares of degrees.
Journal of Graph Theory 31:283--295 (1999)

\bibitem{PS}
Peled, U.N., Srinivasan, M.K.:
The polytope of degree sequences.
Linear Algebra and its Applications 114/115:349--377 (1989)

\bibitem{Rys}
Ryser, H.J.:
Combinatorial properties of matrices of zeroes and ones.
Canadian Journal of Mathematics 9:371--377 (1957)

\bibitem{Sch}
Schrijver A.:
Combinatorial Optimization. Springer (2003)

\bibitem{Wes}
West, D.B.: Introduction to graph theory. Prentice Hall (1996)

\end{thebibliography}
\end{document}